\documentclass[12pt]{article}

\usepackage{amsmath}
\usepackage{amscd}
\usepackage{amsfonts}
\usepackage{epsfig}
\usepackage{theorem}
\usepackage{authblk}
\usepackage[titletoc]{appendix}
\usepackage[margin=1.0in]{geometry}

\def \C {\mathbb C}

\def \R {\mathbb R}

\def \N {\mathbb N}

\def \Z{\mathbb Z}

\title{Properties of Generalized Bessel Functions}

\author{Parker Kuklinski, David A. Hague}

\begin{document}


\maketitle

\begin{abstract}
The Generalized Bessel Function (GBF) extends the single variable Bessel function to several dimensions and indices in a nontrivial manner. Two-dimensional GBFs have been studied extensively in the literature and have found application in laser physics, crystallography, and electromagnetics. In this article, we document several properties of $m$-dimensional GBFs including an underlying partial differential equation structure, asymptotics for simultaneously large order and argument, and analysis of generalized Neumann, Kapteyn, and Schl\"{o}milch series. We extend these results to mixed-type GBFs where appropriate.
\end{abstract}

\section{Introduction}

Bessel functions \cite{watson22} are pervasive in mathematics and physics and are particularly important in the study of wave propagation. Bessel functions were first studied in the context of solutions to a second order differential equation known as Bessel's equation:
\begin{equation}
x^2f''(x)+xf'(x)+(x^2-n^2)f(x)=0.
\label{eq:besselEquation}
\end{equation}
This equation is parameterized by the value $n$.  Solutions to this equation are known as Bessel functions of order $n$. Since equation \eqref{eq:besselEquation} is a second order linear differential equation, there exist two linearly independent solutions. These solutions, denoted $J_n (x)$ and $Y_n (x)$, are referred to as Bessel functions of the first and second kind respectively, where $J_n (0)$ is finite and $Y_n (x)$ has a singularity at $x=0$. These functions commonly arise in physical problems involving cylindrical equations including the Laplace, Helmholtz, and Schrodinger equations \cite{farlow93}.

For integer order $n$, Bessel functions of the first kind admit an integral representation:
\begin{equation}
J_n(x)=\frac{1}{2\pi}\int _{-\pi}^\pi e^{i(n\theta -x\sin\theta )}d\theta
\label{eq:bfInt}
\end{equation}
Recognizing equation \eqref{eq:bfInt} as the coefficients to a complex Fourier series expansion, we recover the following equation:
\begin{equation}
e^{-ix\sin{\theta}}=\sum _{n=-\infty}^\infty J_n(x)e^{-in\theta}
\label{eq:jacobiAnger}
\end{equation}
Equation \eqref{eq:jacobiAnger} is referred to as a Jacobi-Anger expansion \cite{erdelyi53} and essentially states that the complex Fourier series coefficients for the function $e^{-ix\sin\theta}$ are in fact the $n^{th}$ order cylindrical Bessel functions of the first kind.  The Bessel functions are oscillatory functions of the variable $x$ with even and odd symmetry in $x$ for even and odd orders $n$ respectively.

Now consider a complex exponential function whose argument itself is represented as a Fourier sine series \cite{HagueDiss, HagueI}, represented as:
\begin{equation}
s(\theta )=\exp \left[ -i\sum_{k=1}^m x_k\sin{k\theta}\right]
\end{equation}
Further suppose we wish to find the complex Fourier series coefficients for that function.
By leveraging a more general form of the Jacobi-Anger expansion, we may write the following representation of $s(\theta )$ \cite{dattoli96} as
\begin{equation}
s(\theta )=\sum _{n=-\infty}^\infty J_n^{\bf p}({\bf x})e^{-in\theta}
\end{equation}
where ${\bf x}:\{ 1,...,m\}\rightarrow\R$, ${\bf p}:\{ 1,...,m\}\rightarrow\N$ with ${\bf p}=\{ p_1,...,p_m\}$, and
\begin{equation}
J_n^{\bf p}({\bf x})=\frac{1}{2\pi}\int _{-\pi}^\pi\exp \left[ i\left( n\theta -\sum_{k=1}^m x_k\sin{p_k\theta}\right)\right] d\theta .
\label{eq:equation6}
\end{equation}
We call $J_n^{p1,\dots, p_m}({\bf x}) = J_n^{\bf p}({\bf x})$ in \eqref{eq:equation6} an $m$-dimensional Generalized Bessel Function (GBF) with indices ${\bf p}$ and require that $\text{gcd}({\bf p})=1$ to avoid trivial simplifications. If any of the arguments in ${\bf x}$ are set to zero, then \eqref{eq:equation6} simplifies to a lower dimensional GBF. We can extend this definition to the mixed-type generalized Bessel functions (MT-GBF) such that the argument also includes cosines \cite{dattoli96}:
\begin{equation}
J_n({\bf x};{\bf y})= \frac{1}{2\pi}\int _{-\pi}^\pi\exp \left[ i\left( n\theta -\sum _{k=1}^m (x_k\sin{k\theta}+y_k\cos{k\theta})\right)\right] d\theta .
\end{equation}
An infinite-dimensional variant of the MT-GBF (and the GBF) can be constructed by requiring ${\bf x},{\bf y}\in\ell ^2(\N )$ such that the integral converges \cite{dattoli1997}.

\begin{figure}[h]
\includegraphics[width=1.0\textwidth]{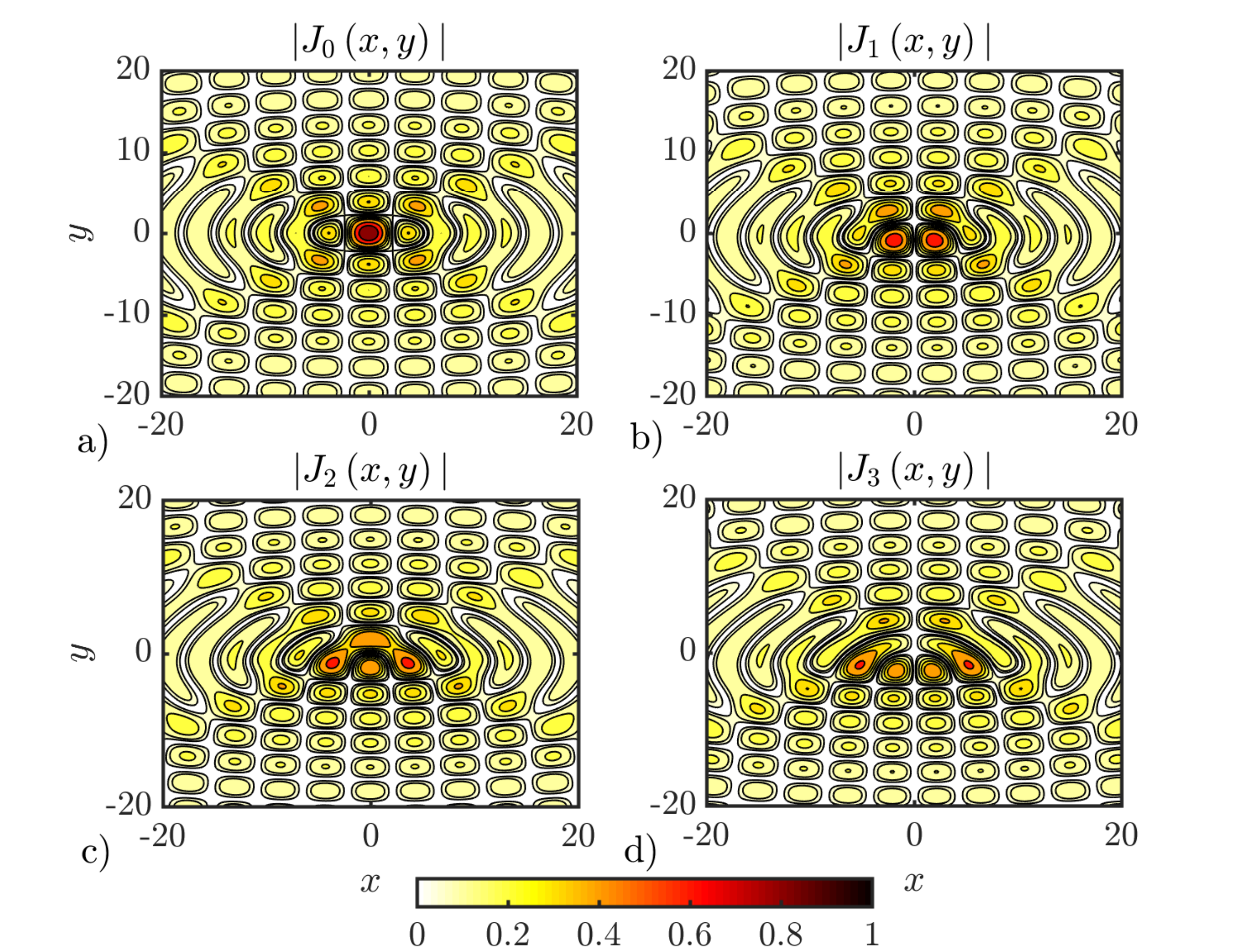}
\caption{Plots of $|J_n^{1,2}(x,y)|$ for $-20<x,y<20$ for (a) $n=0$, (b) $n=1$, (c) $n=2$ and (d) $n=3$. Much like their one-dimensional counterparts, the GBFs are oscillatory in both $x$ and $y$ dimensions. They additionally display symmetries over the $x$ and $y$ axis.}
\label{fig:GBF2}
\end{figure}

\begin{figure}[h]
\includegraphics[width=1.0\textwidth]{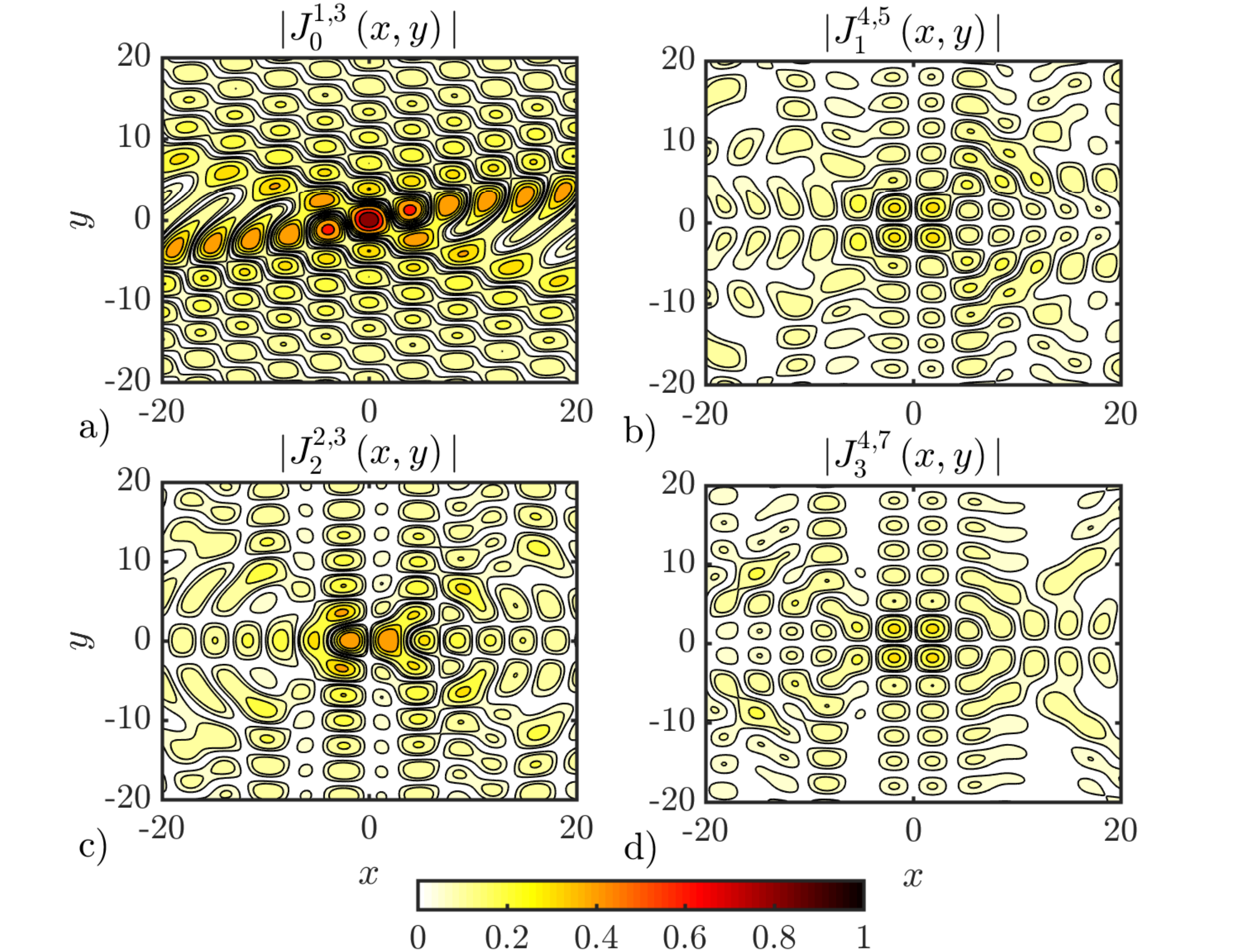}
\caption{Plots of $|J_n^{p,q}(x,y)|$ for $-20<x,y<20$ for (a) $n=0$, $p=1$, $q=3$, (b)$n=1$. $p=4$, $q=5$, (c) $n=2$, $p=2$, $q=3$, (d) $n=3$, $p=4$, $q=7$. These GBFs with their more general indices $p$ and $q$ sill possess similar oscillatory properties as their $J_n^{1,2}$ counterparts with more complex symmetry and structure.}
\label{fig:GBF3}
\end{figure}

GBFs first appeared in the literature in the early 1900s as a simple but isolated extension of the Bessel functions \cite{appell15}.  The 2-D GBFs were extensively utilized in the physics community \cite{Reiss_1962, nikishov1964quantum, Reiss_1980, Reiss_1983} and their mathematical properties were thoroughly explored by the efforts of \cite{dattoli1990theory, dattoli1991theory, dattoli1992linear}.  These functions were re-examined and further generalized to the m-dimensional and infinite-dimensional \cite{dattoli1995, dattoli96, dattoli1997} cases as more of an effort was made to connect the GBF to other notable special functions \cite{dattoli1992generating, andrews85}.  Additionally, the authors have taken an interest in the GBFs as they provide insight into the analysis and design of frequency modulated signals for radar and sonar applications \cite{HagueDiss, HagueI, Hague_MTFFM, hague2020adaptive}.  However, many of the well understood properties and identities of $1-$dimensional Bessel functions have not been extended to the $m-$dimensional GBFs.  In this paper we concentrate on three fundamental properties of GBFs that to the best of the authors' knowledge have not been fully established.  These are (1) finding a second order linear Partial Differential Equation (PDE) of which the GBF is a solution, (2) deriving asymptotic approximations of GBFs with arbitrary dimension for large arguments and/or orders with general indicies ${\bf p}$, and (3) establishing analytical expressions for infinite summations of GBFs such as the Neumann, Kapteyn and Schl\"{o}milch series.  These three fundamental properties for 1-D Bessel functions have found extensive use in Physics.  It is the hope of the authors that extending these properties for m-dimensional GBFs could potentially provide further insight into problems of interest to the Physics community at large.

While the Bessel function was defined as a solution to a differential equation, it is more difficult to discern whether the GBF satisfies a similar partial differential equation. Previous efforts \cite{dattoli96} have showed that the $J_n^{1,m}$ GBF satisfies an order $m^2$ linear PDE which does not bear any passing similarity to Bessel's differential equation.  Work by the authors in \cite{kuklinski18} determined that the GBF does not solve a second order linear PDE but certain modifications of GBFs do; in particular the MT-GBF satisfies a first order linear Schrodinger type equation.  A similar Schrodinger type equation result was obtained by \cite{cesarano} for the special case of 2-dimensional 1-parameter GBFs.  However, this PDE only describes the simple cross sectional structure between two variables $x_k$ and $y_k$ and ignores the more intricate interplay between just the ${\bf x}$ variables.  In this paper, we demonstrate that it is indeed possible to write a second order linear PDE for $J_n^{1,2}\left(x,y\right)$ of which Bessel's original linear PDE is a special case.  

Most analysis of GBFs has been limited to the two dimensional case, and for a complete asymptotic development, even further limited to the case of $(p_1,p_2)=(1,2)$ \cite{korsch06}. Using techniques from algebraic geometry \cite{cox06}, we describe a general method to analytically compute bifurcation surfaces for both GBFs with large argument and for GBFs with simultaneously large order and argument. These bifurcation surfaces are revealed to be piecewise algebraic.  We conclude this paper with a discussion of various summations involving MT-GBFs. In \cite{watson22}, Neumann, Kapteyn, and Schl\"{o}milch series were developed for one-dimensional Bessel functions. We consider generalizations of these series where the MT-GBF replaces the Bessel function. It is possible to write closed forms for the moments of the Neumann and Kapteyn series, while for the Schl\"{o}milch series we can compute smoothness boundaries. Both the region of convergence for the Kapteyn series and the smoothness boundaries of the Schl\"{o}milch series are connected in an algebraic sense to the region of polynomial decay for MT-GBFs of increasing order and argument.

The rest of the paper is organized as follows; in Section \ref{sec:PDE} we write linearly independent PDE representations for $J_n^{1,2}(x,y)$ and note a possible application to computing level sets. In Section 3 we discuss asymptotics of the GBF and provide closed form results. We then develop closed forms for generalized Neumann and Kapteyn series in Section 4, and also describe smoothness boundaries for the generalized Schl\"{o}milch series.  Finally, Section 5 concludes the paper.

\section{Partial Differential Equations Representation}
\label{sec:PDE}

In this section, we demonstrate that the index (1,2) GBF $J_n^{1,2}(x,y)$ solves a set of linearly independent partial differential equations. These equations can be derived from a combination of trigonometric identities and an integral identitiy. 

\subsection{Partial Differential Equations for the 2-D GBF and m-D MT-GBF}
If $f:\R\rightarrow\R$ is a differentiable periodic function with period $P$, then
\begin{equation}
\int _{x_0}^{x_0+P}f'(x)dx=0.
\label{eq:equation8}
\end{equation}
Let $h(\theta )=n\theta -x\sin\theta -y\sin 2\theta$, $f(x,y)=2\pi J_n^{1,2}(x,y)=\int _{-\pi}^\pi e^{ih(\theta )}d\theta$, and $g(x,y)=\int _{-\pi}^\pi \cos\theta e^{ih(\theta )}d\theta$. To take derivatives of these expressions, we can interchange the derivative and integral operators; computing these derivatives involves integrating products of trigonometric functions and $e^{ih(\theta )}$. This allows us to exploit trigonometric identities to relate the derivatives of $f$ and $g$; for instance by $\sin 2\theta=2\sin\theta\cos\theta$ we have $f_y=2g_x$.

We apply equation \eqref{eq:equation8} to $e^{ih(\theta )}$, $(\sin\theta )e^{ih(\theta )}$, and $(\cos\theta )e^{ih(\theta )}$ to generate the partial differential equations. The first identity becomes
\begin{equation}
0=nf-xg-2y\int _{-\pi}^\pi (\cos 2\theta )e^{ih(\theta )}d\theta
\end{equation}
Using the identity $\cos 2\theta =1-2\sin^2\theta$, we have
\begin{equation}
0=(n-2y)f-xg-4yf_{xx}
\label{eq:equation10}
\end{equation}
This allows us to represent $g$ in terms of $f$ and $f_{xx}$.  Applying equation \eqref{eq:equation8} to $(\sin\theta )e^{ih(\theta )}$ leads to the following formula:
\begin{equation}
0=g-nf_x+\frac{1}{2}xf_y-2iy\int_{-\pi}^\pi(\sin\theta \cos 2\theta )e^{ih(\theta )}d\theta
\end{equation}
Using the identity $\sin\theta\cos 2\theta =\sin 2\theta\cos\theta -\sin\theta$, we can write:
\begin{equation}
0=g-(n+2y)f_x+\frac{1}{2}xf_y+2yg_y
\end{equation}
By taking a derivative of this equation with respect to $x$ and using the identities $f_y=2g_x$ and $f_{yy}=2g_{xy}$, we can write the first partial differential equation:
\begin{equation}
(n+2y)f_{xx}-yf_{yy}-\frac{x}{2}f_{xy}-f_y=0
\label{eq:eq2_6}
\end{equation}

We can derive the second partial differential equation by applying equation \eqref{eq:equation8} to $(\cos\theta )e^{ih(\theta )}$. This leads us to the following formula:
\begin{equation}
0=f_x-ng+x\int _{-\pi}^\pi (\cos ^2\theta )e^{ih(\theta )}d\theta + 2y\int _{-\pi}^\pi (\cos\theta\cos 2\theta )e^{ih(\theta )}d\theta
\end{equation}
By using trigonometric substitutions, we can write these integrals in terms of the derivatives of $f$ and $g$:
\begin{equation}
0=xf_{xx}+f_x+xf-(n-2y)g+2yf_{xy}
\end{equation}
Using the previous substitution for $g$ in equation \eqref{eq:equation10}, we have:
\begin{equation}
0=[x^2+4y(n-2y)]f_{xx}+2xyf_{xy}+ xf_x+[x^2-(n-2y)^2]f
\label{eq:diffEQ}
\end{equation}
There are several promising similarities between equation \eqref{eq:diffEQ} and Bessel's differential equation in equation \eqref{eq:besselEquation}, particularly that the correct degree polynomial in $(x,y)$ multiplies the corresponding ordered derivative, and the dependence on $n^2$ in the coefficient of $f$. Moreover, substituting $y=0$ returns Bessel's differential equation.

We can use a similar procedure to find linearly independent PDEs which govern MT-GBFs of arbitrary dimension. Note that by $\sin^2\theta +\cos ^2\theta =1$, we have
\begin{equation}
\frac{\partial ^2J_n}{\partial x_k^2}+\frac{\partial ^2J_n}{\partial y_k^2}+J_n=0
\end{equation}
This shows that every $(x_k,y_k)$ plane has circular level sets when all other variables fixed. It was shown in \cite{kuklinski18, cesarano} that the MT-GBF also satisfies a Schrodinger-type PDE:
\begin{equation}
nJ_n=i\sum _{k=1}^m k\left( x_k\frac{\partial J_n}{\partial y_k}-y_k\frac{\partial J_n}{\partial x_k}\right)
\end{equation}
A complete list of linearly independent second order PDE identities for the $J_n^{1,2}$ MT-GBF could be compiled using equation \eqref{eq:equation8}, however since there are four input variables we omit this list for brevity.

\subsection{A Comment on Level Sets}

If we could find another second order linear PDE which $J_n^{1,2}(x,y)$ solves independent of \eqref{eq:eq2_6} and \eqref{eq:diffEQ}, then it would be possible to parameterize its vanishing level sets using a second order variant of the method of characteristics \cite{evans98}. Suppose we have smooth functions $(x(t),y(t))$ such that $f(x(t),y(t))=J_n^{1,2}(x(t),y(t))=0$. Then $x'(t)f_x+y'(t)f_y=0$, and taking a second derivative of this equation gives
\begin{equation}
(x')^2f_{xx}+2x'y'f_{xy}+(y')^2f_{yy}+x''f_x+y''f_y=0
\end{equation}
This allows us to write the following matrix representation of the system:
\begin{equation}
\begin{bmatrix} 2(n+2y) & -x & -2y & 0 & -2 \\ x^2+4y(n-2y) & 2xy & 0 & x & 0 \\ (x')^2 & 2x'y' & (y')^2 & x'' & y'' \\ 0 & 0 & 0 & x' & y'\end{bmatrix}\begin{bmatrix} f_{xx} \\ f_{xy} \\ f_{yy} \\ f_x \\ f_y\end{bmatrix}=\begin{bmatrix} 0 \\ 0 \\ 0 \\ 0\end{bmatrix}
\end{equation}
We could then write a system of second order nonlinear ordinary differential equations with initial conditions $(x(0),y(0),x'(0),y'(0))=(j_{n,k},0,0,1)$ or for even $n$ $(x(0),y(0),x'(0),y'(0))=(0,j_{n/2,k},1,0)$ where $j_{n,k}$ is the $k^\text{th}$ root of the one-dimensional Bessel function $J_n(x)$.  It may be possible to determine the topology of the level sets of $J_n^{1,2}(x,y)$ by analyzing this ODE system. Note in Figure \ref{fig:GBF2} that for even $n$, the zero surfaces intersecting the $y$-axis appear to be closed loops, while for $n$ odd there appears to be a single infinite contour winding about the $y$-axis.  This is a potential topic for future investigation.

\section{Asymptotic Properties}

In this section, we characterize the bifurcation surfaces of the GBF using the method of stationary phase \cite{stein93}. Consider the integral:
\begin{equation}
I(t)=\int_a^b g(\theta )e^{itf(\theta )}d\theta
\label{eq:equation21}
\end{equation}
where $g$ and $f$ are smooth functions compactly supported in the interval. Consider the set $S=\{ x\in (a,b):f'(x)=0,f''(x)\ne 0\}$ whose members we refer to as \emph{points of stationary phase}. If this set is nonempty, then we can make an approximation on $I(t)$ for large $t$:
\begin{equation}
I(t)=\sum _{x\in S}g(x)e^{itf(x)}\sqrt{\frac{\pi}{t|f''(x)|}}(1\pm i)+O(t^{-1})
\end{equation}
Here, the phase of the summand (i.e. the sign of $1\pm i$) is determined by the sign of $f''(x)$. If the set of stationary phase points is empty, then $I(t)$ decays superpolynomially, and under minor conditions, exponentially. Some care must be taken that there are no stationary phase points at the endpoints of the integral, and that the stationary phase points do not give the phase function a vanishing second derivative. The latter type of stationary phase points are referred to as \emph{critical points}, and crossing over these points changes the structure of the approximation (i.e. changing the number of stationary phase points in the region of integration).

Compare the integral in equation \eqref{eq:equation21} with our definition of the GBF in equation \eqref{eq:equation6}. For an $m$-dimensional GBF with fixed indices, there are $m+1$ terms which can vary (the order and $m$ input arguments), any combination of which we can choose to be large and apply a stationary phase approximation to. In this section we will consider two cases; we first let the arguments be large relative to the order, and then we let the arguments and the order be simultaneously large. In practice, when the indices are arbitrary integers, these approximations can be difficult to analytically display as they involve solving high order polynomials. However, we can analytically find the locus of critical points which trace out \emph{bifurcation surfaces} in $m$-dimensional space. For relatively small order these bifurcation surfaces become two linear $m-1$ dimensional surfaces, but for large order the bifurcation surfaces include an additional algebraic surface. The collection of bifurcation surfaces of the MT-GBF generally contains only higher order algebraic surfaces.

These more complex algebraic surfaces satisfy systems of polynomial equations which can be solved using the Sylvester matrix determinant. If we have a system of polynomial equations $f(x)=0$ and $g(x)=0$ with coefficients in the set $\{ a_1,...,a_n\}$ and which have no common factors, then there exists a nontrivial multivariate polynomial called the \emph{resultant} which satisfies $\text{Res}(f,g;x)=F(a_1,...,a_n )=0$ \cite{cox06}. If $f(x)=\sum _{k=0}^n a_kx^k$ and $g(x)=\sum _{k=0}^n b_kx^k$, the resultant may be written as the determinant of the Sylvester matrix:
\begin{equation}
\text{Res}(f,g;x)=\det\begin{bmatrix} a_n & a_{n-1} & \hdots & a_0 & 0 & ~ & ~ & ~ \\ b_n & b_{n-1} & \hdots & b_0 & 0 & \ddots & ~ & ~ \\ 0 & a_n & \hdots & a_1 & a_0 & \ddots & ~ & ~ \\ 0 & b_n & \hdots & b_1 & b_0 & \ddots & ~ & ~ \\ ~ & ~ & \ddots & \ddots & \ddots & \ddots & ~ & ~ \\ ~ & ~ & ~ & ~ & ~ & b_n & \hdots & b_0\end{bmatrix}
\end{equation}
Suppose we have a system of $m+1$ polynomial equations in $m$ variables such that $f_{1,k}(x_1,...,x_m)=0$ for $k\le m+1$. Then we can generate a system of $m$ polynomial equations $f_{2,k}(x_1,...,x_{m-1})=\text{Res}(f_{1,k},f_{1,m+1};x_m)$ and continue iteratively until there is one polynomial with none of the $\{ x_k\}$ variables. This polynomial is the resultant of the entire system.

One caveat of this method is that the equation $\text{Res}(f,g)=F(a_1,...,a_n)=0$ will include solutions $\{a_1,...,a_n\}$ which do not simultaneously satisfy $f(x)=0$ and $g(x)=0$ for some $x$, whereas if the system is satisfied for some $x$, the coefficients must solve $F(a_1,...,a_n )=0$. For example, consider the system $f(x)=ax+b+1=0$ and $g(x)=cx+d=0$. For this system to be simultaneously satisfied for some $x$, we must have $F(a,b,c,d)=ad-bc-c=0$. The trivial solution $\{ a,b,c,d\} =\{ 0,0,0,0\}$ satisfies this equation but clearly does not satisfy the system $f\left(x\right) = g\left(x\right) = 0$ for any $x$. Moreover, it will often occur that we would like for the simultaneous solution $x$ to be real or otherwise satisfy some chosen constraint. Only a portion of the bifurcation curve will correspond to this situation, and this information is not encoded in the resultant. We will examine these cases as they arise.  We also note that while we can solve the bifurcation surfaces analytically, these do not give a complete asymptotic description of the GBF. This would require us to solve for the points of stationary phase which satisfy a polynomial of degree $\text{max }{\bf p}$ in $\cos\theta$; which is typically not analytically possible.

\subsection{Large Arguments}
\label{subsec:largeArgs}

\begin{figure}[!h]
\includegraphics[width=1.0\textwidth]{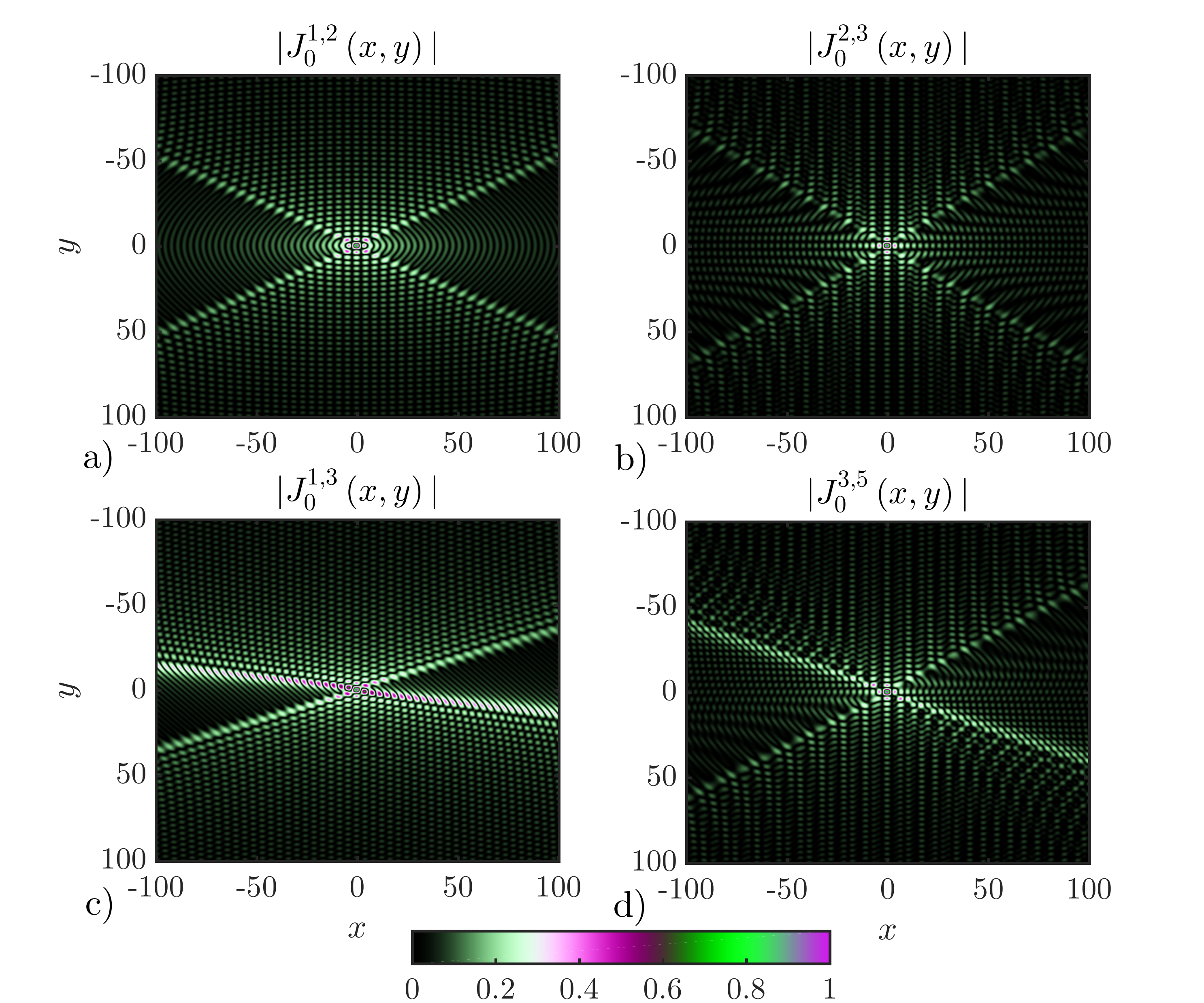}
\caption{Plots of $|J_0^{p,q}(x,y)|$ for $-20<x,y<20$ for (a) $p=1$, $q=2$, (b) $p=2$, $q=3$, (c) $p=1$, $q=3$ and (d) $p=3$, $q=5$. Notice that the top two figures, each containing an even index, have bifurcation curves which are symmetric about both coordinate axes. Meanwhile, the bottom two figures, both with odd indices, contain bifurcation curves with only rotational symmetry.}
\label{fig:GBF4}
\end{figure}

We first consider bifurcation curves of the GBF with large arguments relative to the order, as first elucidated by Korsch et. al. \cite{korsch06}. Otherwise stated, let us consider the function $J_n^{{\bf p}} (tx_1,...,tx_m )$ for large values of $t$. We write this function in a stationary phase form
\begin{equation}
J_n^{{\bf p}}(tx_1,...,tx_m)=\frac{1}{2\pi}\int _{-\pi}^\pi e^{in\theta}\exp\left[ -it\sum _{k=1}^m x_k\sin{p_k\theta}\right] d\theta
\label{eq:equation25}
\end{equation}
where $g(\theta )=e^{in\theta}$ and $f(\theta )=-\sum _{k=1}^m x_k\sin{p_k\theta}$. If $(x_1,...x_m)$ is an element of the bifurcation surface, there must exist some $\theta\in [-\pi ,\pi ]$ such that $f'(\theta )=f''(\theta )=0$. Otherwise,
\begin{align}
f'(\theta )&=\sum _{k=1}^m x_kp_k\cos{p_k\theta}=0, \\ f''(\theta )&=\sum _{k=1}^m x_kp_k^2\sin{p_k\theta}=0.
\label{eq:equation26}
\end{align}
We note that these functions can be written as polynomials in $\cos \theta$ and $\sin \theta$. In particular, the Chebyshev polynomials of the first and second kinds respectively satisfy the identities
\begin{align}
T_n(\cos\theta )&=\cos{n\theta}, \\ U_{n-1}(\cos\theta )&=\frac{\sin{n\theta}}{\sin\theta}.
\label{eq:equation27}
\end{align}
Using these identities, we can rewrite the system of equations \eqref{eq:equation26}:
\begin{align}
\sum _{k=1}^m x_kp_k T_{p_k}(\cos\theta )&=0, \\ \sin\theta\left[\sum _{k=1}^m x_kp_k^2U_{p_k-1}(\cos\theta )\right] &=0.
\label{eq:equation28}
\end{align}
The right equation in \eqref{eq:equation28} is trivially satisfied if $\sin\theta =0$, or when $\theta =0$ or $\theta =\pi$. These two solutions lead to the following representations of bifurcation surfaces when plugged back into the left equation of \eqref{eq:equation28}:
\begin{align}
\sum _{k=1}^m x_kp_k&=0, \\ \sum _{k=1}^m x_kp_k(-1)^{p_k}&=0.
\label{eq:equation29}
\end{align}
Notice that if all of the $p_k$ are odd (they cannot all be even or else their greatest common denominator would be at least 2), then both equations of \eqref{eq:equation29} denote the same surface. In this case, letting $\theta =\frac{\pi}{2}$ or $\theta =\frac{3\pi}{2}$ trivially satisfies the left equation of \eqref{eq:equation28}. Substituting these values into the right equation of \eqref{eq:equation28} leads to the following alternate surface:
\begin{equation}
\sum _{k=1}^m x_kp_k^2 (-1)^\frac{p_k+1}{2}=0.
\end{equation}
We illustrate this dichotomy in behavior in Figure \ref{fig:GBF4}; notice that the bifurcation curves of $J_n^{1,2}(x,y)$ and $J_n^{2,3}(x,y)$ are symmetric about the coordinate axes, but the bifurcation curves of $J_n^{1,3}(x,y)$ and $J_n^{3,5}(x,y)$ have only rotational symmetry about the origin.  For the MT-GBF, there is no factor of $\sin\theta$ in $f''(\theta )$, so the bifurcation surfaces will be higher order algebraic surfaces. The bifurcation surface of the $J_n^{1,2}(tx_1,tx_2;ty_1,ty_2)$ MT-GBF is an eighth order algebraic surface whose equation is too large for the page.

\subsection{Large Order and Arguments}
\label{subsec:largeOrdandArgs}
We now compute bifurcation curves of the GBF with simultaneously large order and large arguments, otherwise $J_{tn}^{\bf p} (t{\bf x})$ for large values of $t$. In these cases, there is a region of exponential decay containing the origin which grows linearly with $n$. The bifurcation surfaces will bound this region.  In the stationary phase representation, we have $f(\theta )=n\theta -\sum_{k=1}^m x_k\sin p_k\theta$ such that the bifurcation curves now satisfy:
\begin{align}
f'(\theta )&=\sum _{k=1}^m x_kp_k\cos{p_k\theta}=n, \\ f''(\theta )&=\sum _{k=1}^m x_kp_k^2\sin{p_k\theta}=0.
\label{eq:equation30}
\end{align}
As stated previously, we will be able to write a multivariate polynomial equation in $\{ x_k\}$, a subset of whose solutions satisfy equations \eqref{eq:equation30} for some $\theta$. By factoring out $\sin\theta$ from the second equation, we are able to derive two linear bifurcation surfaces:
\begin{align}
\sum _{k=1}^m x_kp_k&=n, \\ \sum _{k=1}^m x_kp_k(-1)^{p_k}&=n.
\end{align}
We refer to these solutions as the \emph{trivial solutions}. Notice that unlike the small order case of the previous section, these two surfaces are distinct even if ${\bf p}$ contains only odd indices, in which case they are parallel. If ${\bf p}=\{ {\bf p_e},{\bf p_o}\}$ corresponds to the even and odd indices respectively, then the intersection of the trivial surfaces can be represented as the intersection of two orthogonal surfaces:
\begin{align}
\sum _{k\in {\bf p_e}}x_kp_k&=0, \\ \sum _{k\in {\bf p_o}}x_kp_k&=n.
\end{align}

\begin{figure}[!h]
\includegraphics[width=1.0\textwidth]{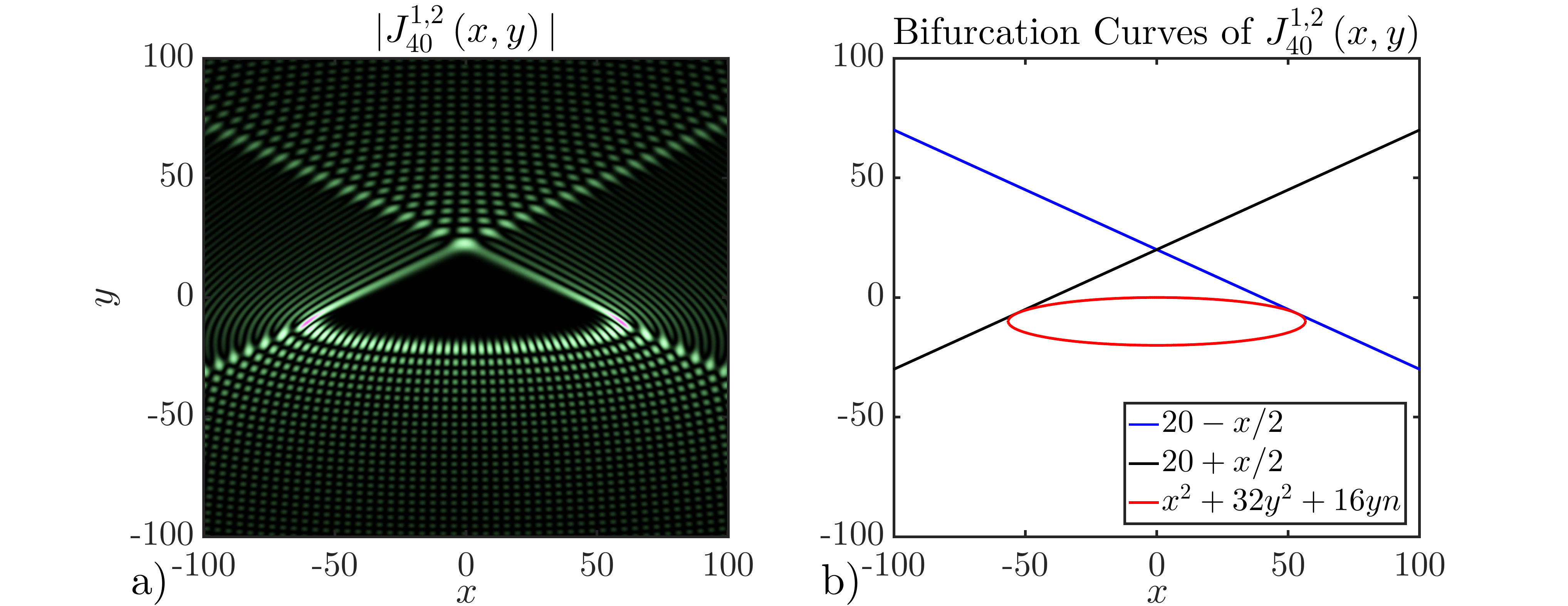}
\caption{Plot of $|J_{40}^{1,2}(x,y)|$ for $-100<x,y<100$ (a) and its corresponding Bifurcation curves (b). Notice that crossing the upper part of the ellipse does not lead to a difference in asymptotic behavior.}
\label{fig:GBF5}
\end{figure}

The nontrivial bifurcation surfaces can be computed using the resultant as described at the beginning of the section. For instance, the nontrivial bifurcation curve of $J_{tn}^{1,2}(tx,ty)$ becomes
\begin{equation}
x^2+32y^2+16yn=0
\label{eq:equation33}
\end{equation}
which is the result obtained by \cite{lotstedt09}.  We plot the bifurcation curves of this GBF in Figure \ref{fig:GBF5}. Notice here that the upper section of the ellipse as predicted by equation \eqref{eq:equation33} does not actually behave as a bifurcation curve, i.e. it subdivides a coherent region of exponential decay. We remedy this discrepancy by noting that the simulateous solutions of equation \eqref{eq:equation30} must satisfy $|\cos\theta |\le 1$, otherwise $\theta$ is not in the region of integration. In this example, this implies that we must have $\left\lvert\frac{x}{8y}\right\rvert\le 1$, or equivalently the section of the ellipse which lies above the ine $y=-\frac{n}{6}$ is not truly part of the bifurcation curve.

Since previous asymptotic analysis of GBFs \cite{korsch06} \cite{lotstedt09} was dependent on explicitly solving for stationary phase points, those authors could only consider the $J_{tn}^{1,2}(tx,ty)$ GBF because the stationary phase points satisfy a quadratic function in $\cos\theta$. However, the resultant allows us to solve for bifurcation curves of MT-GBFs of arbitrary order. For example, the nontrivial bifurcation curve of $J_{tn}^{1,3}(tx,ty)$ becomes
\begin{equation}
(x-9y)^3+81yn^2=0
\end{equation}
This representation overstates the true bifurcation curve, and we can show that only the section of this curve which satisfies $\frac{1}{4}-\frac{x}{36y}\le 1$ acts as a legitimate bifurcation curve, as shown in Figure 6.  We could additionally compute bifurcation surfaces of the MT-GBF in this way, but the equations would be far too large to display concisely.

\begin{figure}[!h]
\includegraphics[width=1.0\textwidth]{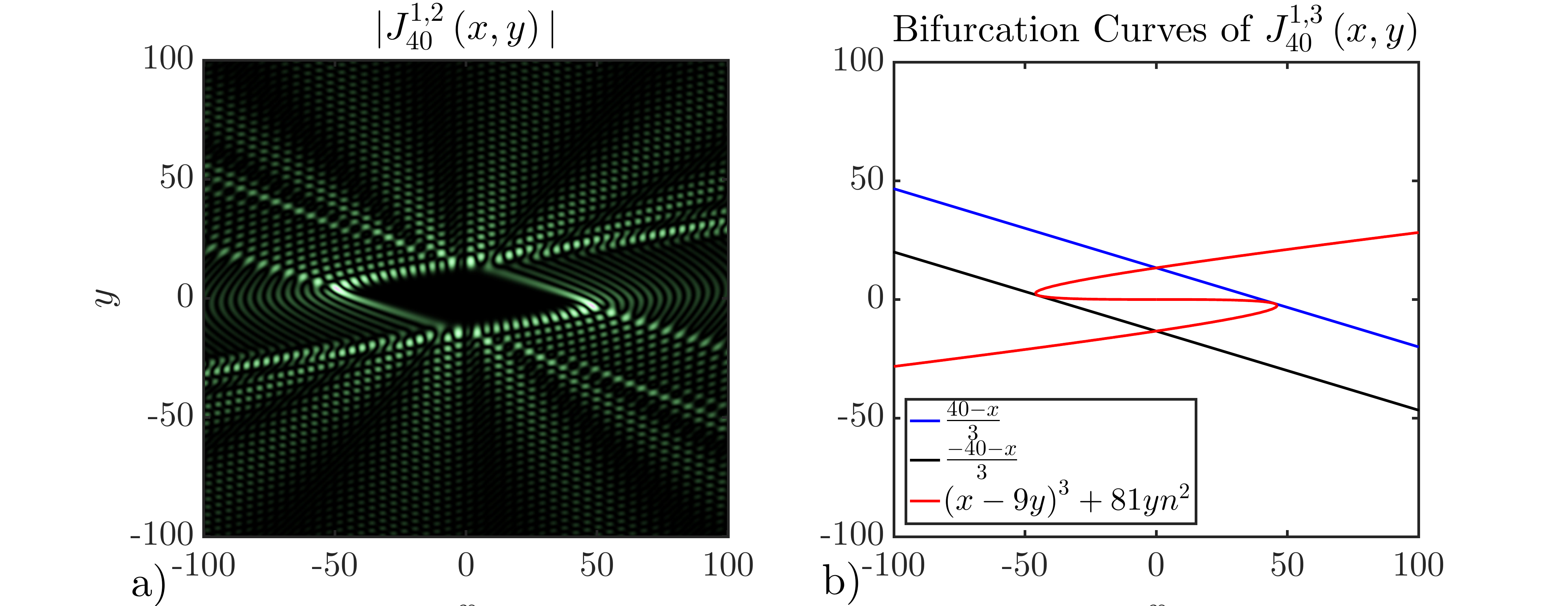}
\caption{Plots of $|J_{40}^{1,3}(x,y)|$ for $-100<x,y<100$ (a) and its corresponding Bifurcation curves (b). Notice that crossing the central part of the cubic curve does not lead to a change in asymptotic behavior.}
\label{fig:GBF6}
\end{figure}

\section{MT-GBF Series}

In this section, we generalize various infinite sums involving the one-dimensional Bessel function to the MT-GBF. The techniques used to compute these sums are not unique to the dimension of the Bessel function and are easily extended to the MT-GBF. For the generalized Neumann and Kapteyn series, we use generating functions to compute $\ell^\text{th}$ moments. These moments could then be used to approximate series involving arbitrary coefficients. The moments of the Schl\"{o}milch series are not so easily calculated even in the one-dimensional case. However, a distinguishing feature of this series in one dimension is that it is only piecewise smooth; the boundaries of these pieces occur at $x=2\pi n$ for $n\in\Z$. We document a similar property in the GBF case and compute the smoothness boundaries using the resultant of the previous section.

\subsection{Generalized Neumann Series}
We first present a compact method to compute the moments of the following MT-GBF based Neumann series
\begin{align}
f_{1,\ell}\left({\bf x};{\bf y}\right) &=\sum _{n=-\infty}^\infty n^{\ell}J_n({\bf x};{\bf y}) \\ 
f_{2,\ell}\left({\bf x};{\bf y}\right) &=\sum _{n=-\infty}^\infty n^{\ell}J_n({\bf x};{\bf y})^2 \\
f_{3,\ell}\left({\bf x};{\bf y}\right) &=\sum _{n=-\infty}^\infty n^{\ell}\left|J_n({\bf x};{\bf y})\right|^2.
\end{align}
These derivations will primarily use the MT-GBF Jacobi-Anger identity as well as the derivative recursion of the MT-GBF
\begin{align}
\dfrac{\partial J_n}{\partial x_k} &= \frac{1}{2}\left(J_{n-k}-J_{n+k}\right)  \\
\dfrac{\partial J_n}{\partial y_k} &= \frac{1}{2i}\left(J_{n-k}+J_{n+k}\right).
\end{align}
We will also utilize the square MT-GBF identity
\begin{align}
J_n\left(x_k;y_k\right)^2 &= \frac{1}{2\pi}\int_{-\pi}^{\pi}J_{2n}\left(2x_k\cos k\theta; 2y_k\cos k\theta\right)d\theta \label{eq:GBF_Square_1} \\ 
\left|J_n\left(x_k;y_k\right)\right|^2 &= \frac{1}{2\pi}\int_{-\pi}^{\pi}J_{2n}\left(2x_k\cos k\theta-2y_k\sin k\theta; 0\right)d\theta. \label{eq:GBF_Square_2}
\end{align}
These integral identities can be readily derived from their 1-D Bessel function counterparts.

\subsubsection{Moments of $f_{1,\ell}\left({\bf x};{\bf y}\right)$}
\label{subsubsec:Series1}
Computing these sums becomes nearly trivial by utilizing the MT-GBF Jacobi-Anger expansion.  Let $g\left(\theta\right) = \sum_{n=-\infty}^{\infty}J_n\left(x_k;y_k\right)e^{-i n\theta}$.  If we let $X_{\ell}=\sum _{n=1}^\infty n^{\ell}x_n$ be the ${\ell}^\text{th}$ moment of ${\bf x}$ and likewise for $Y_{\ell}$, then we can simply represent the first few moments as:
\begin{align}
\mu _0&=e^{-iY_0}, \\ \mu _1&=X_1e^{-iY_0}, \\ \mu _2&=(X_1^2-iY_2)e^{-iY_0}, \\ \mu _3&=(X_3-3iX_1Y_2+X_1^3)e^{-iY_0}
\end{align}
These results are MT-GBF generalizations of the results obtained in \cite{dattoli96}

\subsubsection{Moments of $f_{2,\ell}\left({\bf x};{\bf y}\right)$}
\label{subsubsec:Series2}
For the series $f_{2,\ell}\left({\bf x};{\bf y}\right)$, using the identity in \eqref{eq:GBF_Square_1}, we can interchange the summation and the integral to write
\begin{equation}
f_{2,\ell} = \frac{1}{2\pi}\int_{-\pi}^{\pi}\sum_{-\infty}^{\infty}n^{\ell} J_{2n}\left(2x_k \cos k\theta; 2y_k\cos k\theta\right)d\theta
\label{eq:series2Eq1}
\end{equation} 
Consider the following summation formula
\begin{equation}
h_{\ell}\left(\theta\right) = \sum_{n=-\infty}^{\infty}n^{\ell}J_{2n}\left(x_k;y_k\right)\cos^2n.\theta
\end{equation}
Expanding the $\cos^2n\theta$ in terms of complex exponentials results in the expression
\begin{equation}
h_{\ell}\left(\theta\right) = i^{\ell}\left(\frac{1}{4}f^{\left(\ell\right)}\left(2\theta\right) + \frac{1}{4}f^{\left(\ell\right)}\left(-2\theta\right) + \frac{1}{2}f^{\left(\ell\right)}\left(0\right)  \right)
\end{equation}
recognizing that the following identity holds:
\begin{equation}
 \sum_{n=-\infty}^{\infty}n^{\ell}J_{2n}\left(x_k;y_k\right) = \frac{1}{2^{\ell}}h_{\ell}\left(\pi\right)
\end{equation}
We may now substitute the above result into \eqref{eq:series2Eq1} and conduct the integral, the result of which will be another MT-GBF expression.  Using the fact that $J_n\left(0;y_k\left(-1\right)^k\right) = \left(-1\right)^n J_n\left(0;y_k\right)$, we can derive th following quantities
\begin{align}
f_{2,0} &= J_0\left(0;2y_k\right) \\
f_{2,1} &= i\sum_{k=1}^{\infty}k x_k \dfrac{\partial J_0 \left(0;2y_k\right)}{\partial y_k} \\
f_{2,2} &= \frac{1}{2}\sum_{k=1}^{\infty}k^2y_k\dfrac{\partial J_0\left(0;2y_k\right)}{\partial y_k}- \sum_{j,k=1}^{\infty}jkx_jx_j\frac{\partial^2 J_0\left(0;2y_k\right)}{\partial y_j \partial y_k}
\end{align}

\subsubsection{Moments of $f_{3,\ell}\left({\bf x};{\bf y}\right)$}
\label{subsubsec:Series3}
For the series $f_{3,k}$, again using using the integral identity \cite{eq:GBF_Square_2} we may write
\begin{equation}
f_{3,\ell}=\frac{1}{2\pi}\int_{-\pi}^{\pi}\sum_{n=-\infty}^{\infty}n^{\ell}J_{2n}\left(2x_k\cos k\theta -2y_k\sin k\theta; 0\right).
\end{equation}
Again using the function $h_{\ell}\left(\theta\right)$ from the previous section, we can readily compute the first few moments in this series
\begin{align}
f_{3,0} &= 1 \\
f_{3,1} &= 0 \\
f_{3,2} &= \sum_{k=1}^{\infty}\dfrac{k^2\left(x_k^2+y_k^2\right)}{2}.
\label{eq:series3eq1}
\end{align}
Reassuringly, the result for $f_{3,2}$ in \eqref{eq:series3eq1} collapses back to a well known identity for 1-D Bessel functions for the case where m=1 \cite{watson22},
\begin{equation}
f_{3,2} = \sum_{n=-\infty}^{\infty}n^2|J_n\left(x\right)|^2 = \frac{x^2}{2}.
\end{equation}

\subsection{Generalized Kapteyn Series}
We now compute moments of the generalized Kapteyn series, where the summation index is included both in the order and argument:
\begin{equation}
\mu_{\ell}=\sum _{n=-\infty}^\infty n^{\ell}J_n(n{\bf x};n{\bf y})
\end{equation}
To compute these moments, we develop a generating function from a variant of the multi-dimensional feedback equation \cite{dattoli98} \cite{kuklinski18}. Let $f(\theta )=\theta -\sum _{k=1}^\infty (x_k\sin k\theta +y_k\cos k\theta )$ and suppose we wish to invert this function, i.e. represent $\theta =f^{-1}(t)$ as a function of $t$. Notice that this inversion is only valid on the set $\Omega$ of coordinates $\{ {\bf x},{\bf y}\}$ on which $f(\theta )$ is monotone increasing. It can be shown that $\Omega$ corresponds to the region of exponential decay of the previous section with $n=1$. This is intuitive as for $\{ {\bf x},{\bf y}\}$ outside this region, there will be stationary phase points in the integral represnetation such that $J_n(n{\bf x};n{\bf y})=O(n^{-1/2})$, and therefore the series will not converge.

If we assume $\{ {\bf x},{\bf y}\}\in\Omega$, then both $f$ and $f^{-1}$ are monotone increasing, continuous, bijective, and invertible. Since $f(\theta +2\pi )=f(\theta )+2\pi$, we have
\begin{equation}
f(f^{-1}(t+2\pi ))=t+2\pi =f(f^{-1}(t)+2\pi )
\end{equation}
This implies that $g(t)=f^{-1}(t)-t$ is $2\pi$ periodic such that $g$ admits a Fourier expansion $g(t)=\sum _{m=-\infty}^\infty a_me^{imt}$ and the coefficients satisfy $a_m=\frac{1}{2\pi}\int _{-\pi}^\pi g(t)e^{-imt}dt$. By integration by parts, a change of variables to $\theta$, and recognizing that  \\ $\int _{-\pi}^\pi f'(\theta )e^{-imf(\theta )}d\theta =0$, we can write these coefficients as MT-GBFs:
\begin{equation}
a_{-m}=-\frac{1}{2\pi mi}\int _{-\pi}^\pi e^{imf(\theta )}d\theta =-\frac{i}{m}J_m(m{\bf x};m{\bf y})
\end{equation}

If we define $f(\theta )=\theta -h(\theta )$, then we can write a generating function for the generalized Kapteyn series:
\begin{equation}
h(\theta )=-i\sum _{m\ne 0}\frac{1}{m}J_m(m{\bf x};m{\bf y})e^{im(\theta -h(\theta ))}
\end{equation}
To solve for moments of the Kapteyn series, we must use the quantity $\theta _0$ which satisfies $f(\theta _0)=0$. Since $f$ is bijective, $\theta _0$ exists and is unique. Furthermore if ${\bf y}=0$, then $\theta _0=0$, $h^{(2n)}(\theta _0)=0$, and $h^{(2n+1)}(\theta _0)=\sum k^{2n+1}x_k$. The moments arise by taking derivatives of both sides of equation (40) with respect to $\theta$ and evaluating at $\theta _0$. We list the first few moments here:
\begin{align}
\mu _0&=\frac{h'(\theta _0)}{1-h'(\theta _0)}, \\ \mu _1&=-\frac{ih''(\theta _0)}{(1-h'(\theta _0))^3}, \\ \mu_2&=-\frac{h'''(\theta _0)}{(1-h'(\theta _0))^3}-\frac{3h''(\theta _0)^2}{(1-h'(\theta _0))^4}
\end{align}
These are MT-GBF versions of the equations that appear in Dattoli et. al. \cite{dattoli98}.

\subsection{Generalized Schl\"{o}milch Series}

We conclude this section with a discussion of generalized Schl\"{o}milch series which we write as
\begin{equation}
\mu_{\ell}=\sum _{n=1}^\infty n^{-{\ell}}J_m(n{\bf x};n{\bf y})
\end{equation}
where $m$ is fixed and $\ell$ is a positive integer. In Watson \cite{watson22}, only special cases of the one-dimensional Schl\"{o}milch series admit a simple algebraic representation. However, it is noted that these summations are not smooth in neighborhoods about $x=2\pi n$ for $n\in\Z$. We attempt to recover a similar result in the higher dimensional MT-GBF case.

For $\ell\ge 2$, we can interchange the sum and integral to write
\begin{equation}
\mu _{\ell}=\frac{1}{2\pi}\int _{-\pi}^\pi e^{im\theta}\left[\sum _{n=1}^\infty\frac{e^{-inh(\theta )}}{n^{\ell}}\right] d\theta = \frac{1}{2\pi}\int _{-\pi}^\pi e^{im\theta}\text{Li}_{\ell}(e^{-ih(\theta )})
\end{equation}
where $\text{Li}_{\ell}(z)$ is the polylogarithm function \cite{nielsen09}. The polylogarithm is valid for $|z|<1$ but can be extended to an analytic function in $z \in \C\backslash [1,\infty]$. As such, the integrand of equation (43) passes through the branch point at $z=1$ if there exist $\theta\in [-\pi ,\pi]$ such that $h(\theta )=2\pi n$ for some $n\in\Z$. This creates a collection of surfaces in $({\bf x},{\bf y})$ space which we term \emph{smoothness boundaries}; $\mu _{\ell}$ is not smooth in neighborhoods of the smoothness boundaries and crossing over one of these boundaries changes the number of times the integrand passes through the branch point. Therefore, we index the smoothness boundaries by the simultaneous equations
\begin{equation}
h(\theta )=2\pi n,\hspace{1cm} h'(\theta )=0
\end{equation}
Note that this closely relates to the discussion of bifurcation curves in the previous section. We can use similar methods to compute these smoothness boundaries. For $J_n^{1,2}(x,y)$, these boundaries become:
\begin{equation}
64(2\pi n)^4y^2+(2\pi n)^2(x^4-80x^2y^2-128y^4)+ (64y^6-48x^2y^4+12x^4y^2-x^6)=0
\end{equation}
We display these surfaces in Figure 7 along with the $J_n^{1,3}(x,y)$ smoothness boundaries located at $(2\pi n)^2y-(x^3+9x^2y+27xy^2+27y^3)=0$.

\begin{figure}[!h]
\includegraphics[width=1.0\textwidth]{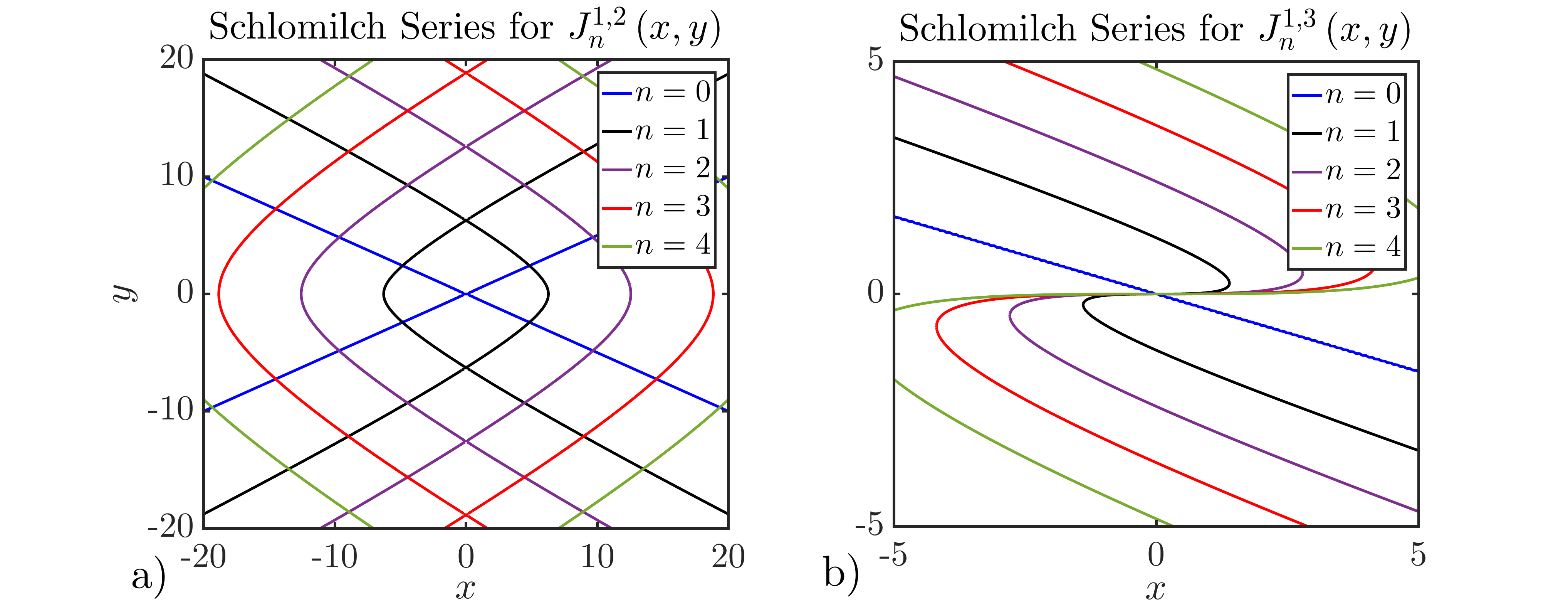}
\caption{(\emph{Left}) Smoothness boundaries of the $J_n^{1,2}(x,y)$  Schl\"{o}milch Series (\emph{Right}) Smoothness boundaries of the $J_n^{1,3}(x,y)$  Schl\"{o}milch series. Moments of the Scl\"{o}milch series are continuous in the domain but not smooth; they are smooth in regions bounded by these curves, however.}
\label{fig:GBF7}
\end{figure}

\section{Conclusion}

In this paper we have documented several novel properties of multi-dimensional GBFs from their one-dimensional counterparts. The properties listed here, while bearing similarity to analogous results on one-dimensional Bessel functions, have a distinct algebraic geometry interpretation. The authors believe that the results in each of these sections can be improved upon with further research. Though it is likely that $J_n^{1,2}(x,y)$ does not solve another linearly independent partial differential equation, it may be possible to solve for the topology of its level sets via the ordinary differential equation system generated from equation (16). It would also be worthwhile to pursue more precise asymptotics for higher dimension GBFs beyond the bifurcation surface calculations presented here.  Since many of these properties of the one-dimensional Bessel functions have found such extensive use in the field of mathematical physics, it is the hope of the authors that these GBF properties might be of similar applicability.


\end{document}